\documentclass[11pt,reqno]{amsart}

\usepackage[utf8]{inputenc}
\usepackage[english]{babel}
\usepackage{amsmath,amsfonts,amssymb,amsthm,shuffle}
\usepackage[T1]{fontenc}
\usepackage{lmodern}
\usepackage{mathtools}
\usepackage[titletoc]{appendix}

\usepackage[top=3.5cm,bottom=3.5cm,left=3.6cm,right=3.6cm]{geometry}

\usepackage[dvipsnames]{xcolor}
\usepackage[hyperindex=true,frenchlinks=true,colorlinks=true,
citecolor=Mahogany,linkcolor=DarkOrchid,urlcolor=Tan,linktocpage,
pagebackref=true]{hyperref}

\usepackage{tikz}
\usetikzlibrary{shapes}
\usetikzlibrary{fit}
\usetikzlibrary{decorations.pathmorphing}

\usepackage{dsfont}
\usepackage{wasysym}
\usepackage{stmaryrd}
\usepackage{cite}
\usepackage{subfigure}
\usepackage{multirow}
\usepackage{enumitem}
\usepackage{multicol}
\usepackage{bm}
\usepackage{amsaddr}
\newcommand{\ds}{\displaystyle}
\title[]{Large deviation inequalities for martingales in Banach spaces}
\date{\today}
\author{Xiequan Fan      and   Davide  Giraudo}
\address{Center for Applied Mathematics,
Tianjin University, Tianjin, China. }
\address{Ruhr-Universit\"{a}t Bochum Fakult\"{a}t f\"{u}r Mathematik, Bochum, Germany.}

\numberwithin{equation}{section}
\renewcommand{\leq}{\leqslant}
\renewcommand{\geq}{\geqslant}

\newtheorem{theorem}{Theorem}[section]

\newtheorem{lemma}[theorem]{Lemma}

\newtheorem{Definition}[theorem]{Definition}
 
\newtheorem{corollary}[theorem]{Corollary}

\theoremstyle{remark}
\newtheorem{remark}[theorem]{Remark}

\tikzstyle{Vertex}=[circle,draw=LimeGreen!80,fill=LimeGreen!8,
inner sep=1pt,minimum size=2mm,line width=1pt,font=\scriptsize]
\tikzstyle{Node}=[Vertex,draw=RoyalBlue!80,fill=RoyalBlue!8,inner sep=1.5pt]
\tikzstyle{Leaf}=[rectangle,draw=Black!70,fill=Black!16,
inner sep=0pt,minimum size=1mm,line width=1.25pt]
\tikzstyle{Edge}=[Maroon!80,cap=round,line width=1pt]
\tikzstyle{Mark1}=[draw=BrickRed!80,fill=BrickRed!8]
\tikzstyle{Mark2}=[draw=BurntOrange!80,fill=BurntOrange!8]
\tikzstyle{EdgeRew}=[->,RedOrange!80,cap=round,thick]

\newcommand{\F}{\mathcal{F}}

\newcommand \ens[1]{\left\{ #1\right\}}

\newcommand \PP{\mathbb P}

\newcommand{\E}[1]{\mathbb E\left[#1\right]}

\newcommand \abs[1]{\left|#1\right|}

\newcommand{\norm}[1]{\left\lVert #1 \right\rVert}

\newcommand{\prt}[1]{\left(#1\right)}

\begin{document}
\begin{abstract}
Let $\prt{X _i, \mathcal{F}_i}_{i\geq1}$ be a martingale difference sequence in a  smooth Banach space.
Let $S_n=\sum_{i=1}^nX_i, n\geq 1,$ be the partial sums of $\prt{X _i, \mathcal{F}_i}_{i\geq1}$.
We give upper bounds on the quantity $\PP\prt{\max_{1\leq k\leq n}\norm{S_k}>nx}$ in terms
of $ n\geq 1$ and $x>0$ in two different situations: when the martingale differences have
uniformly bounded exponential moments and when the decay of the tail of the increments is polynomial.

\end{abstract}

 
\maketitle


\section{Introduction} \label{sec1}

Let $\prt{X_i}_{i\geq 1}$ be a sequence of random variables  with values in a
separable Banach space $\prt{B,\norm{\cdot}}$. Consider the partial sums $S_n:=\sum_{i=1}^nX_i, n\geq 1$.
We are interested in establishing upper bound for the probabilities of large deviation, namely
\begin{equation}
\PP\prt{\max_{1\leq k\leq n}\norm{S_k}>nx}, n\geq 1,x>0.
\end{equation}

In the real valued case, we know that
if $\prt{X_i}_{i\geq 1}$ is centered, strictly stationary and ergodic,
then the quantity $\PP \prt{\max_{1\leq k\leq n}\norm{S_k}>nx}
$ converges to $0$ for all fixed $x$. It is also the case when $\prt{X_i}_{i\geq 1}$ is
 an $\mathbb L^2$ bounded martingale differences sequence.
 In this note, we will generalise to Banach space the known results on the convergence
 rates of the probability of large deviation for real-valued martingale differences sequences, namely:
 \begin{itemize}
 \item in \cite{MR1856684} (see Theorem~3.6 therein), Lesigne and Voln\'{y} have established that if $C_1:=\sup_{i\geq 1} \E{|X_{i}|^{p}}$ is finite for some $p\geq 2$,
 then, for all $ x>0$,
 \begin{equation}
 \PP\prt{\max_{1\leq k\leq n}\abs{S_k}>nx}\leq
 \prt{18 p\sqrt{\frac{p}{p-1}}}^p\frac{M^p}{x^p}n^{-p/2}.
\end{equation}
In particular, it implies
\begin{eqnarray}
\mathbb{P}\left( \max_{1\leq k \leq n} |S_k| >n  \right)  &=&O \left( n^{-p/2 } \right),\ \ \ \ n\rightarrow \infty.
\end{eqnarray}
They also showed that the power $p/2$ in the last equality  is  optimal even for stationary and ergodic martingale difference sequences.

\item in Theorem~2.1 of \cite{MR3005732}, it is proved that if
 $C_2:=\sup_{ i\geq 1 }\E{\exp\left\{ |X_{i}|^{\frac{2\alpha}{1-\alpha}}\right\}  }$ is finite
 for some   $\alpha \in (0, 1)$, then, for all $ x>0$,
\begin{eqnarray}
\mathbb{P}\left( \max_{1\leq k \leq n} S_k >n x \right) &\leq& C(\alpha, x)  \exp\left\{-\left(\frac{ x  }{ 4 }\right)^{2\alpha} n^\alpha\ \right\} ,
\end{eqnarray}
where
\[
C(\alpha,  x)=  2+ 35 C_1  \left( \frac{1}
{ x^{2\alpha} 16^{1-\alpha}}  + \frac{1}{  x^2 }
\left( \frac{3(1-\alpha)}{2\alpha}\right)^{\frac{ 1-\alpha}{\alpha}} \right)
\]
does not depend  on $n$.
In particular, with $x=1$, it implies
\begin{eqnarray}
\mathbb{P}\left( \max_{1\leq k \leq n} S_k >n  \right)  &=&O \left( \exp\left\{ - \frac{1}{16} \, n^{\alpha} \right\} \right),\ \ \ \ n\rightarrow \infty. \label{fnkd}
\end{eqnarray}
It is also showed that the power $\alpha$ in (\ref{fnkd}) is optimal even for stationary martingale difference sequences.  See also \cite{MR1856684} for
the case of $\alpha=1/3.$
 \end{itemize}

The extension of deviation or moment inequalities to Banach spaces
is in general not an easy task. Most of the techniques working in the
real-valued which led to deviation inequalities \cite{MR2021875,MR2462551,MR3311214}
does not seem to extend to Banach space valued martingales. Nevertheless, some inequalities are
available in this context: see \cite{1603.00432} for polynomial decay and
\cite{MR1331198} for exponential inequality. This appears to be the adapted tool
for the control of the large deviation probabilities. Here we will be concerned in
two cases: when the exponential moments of the martingale differences are finite and when the conditional moments of the martingale differences are finite.  

Before we state the results, we need to define the notion of Banach space valued martingale.

\begin{Definition}
 Let $\prt{\Omega,\mathcal F,\PP}$ be a probability space and
let $\prt{B,\norm{\cdot}_{B}}$ be a separable Banach space. For any
$p\geq 1$, we denote by $\mathbb L^p_B$ the space of $B$-valued
random variables such that $\norm{X}_{\mathbb L^p_B}^p=\E{ \norm{X}^p  }$
is finite. Let $\prt{\F_i}_{i\geq 1}$ be a  non-decreasing sequence of
sub-$\sigma$-algebras of $\F$. We say that a sequence of $B$-valued
random variables $\prt{X_i}_{i\geq 1}$ is a martingale difference
sequence with respect to the filtration $\prt{\F_i}_{i\geq 1}$ if
\begin{enumerate}
 \item for any $i\geq 1$, $X_i$ is $\F_i$-measurable and
 belongs to $\mathbb L^1_B$;
 \item for any $i\geq 2$, $\E{X_i\mid \F_{i-1}}=0$ almost
 surely.
\end{enumerate}
\end{Definition}

For the validity of deviation inequalities, some assumptions on the geometry of the
involved Banach space have to be made.

\begin{Definition}
 Following \cite{MR0394135}, we say that a Banach space
 $\prt{B,\norm{\cdot}}$ is $r$-smooth ($1<r\leq 2$)
 if there exists an equivalent norm $\norm{\cdot}'$ such that
 \begin{equation*}
  \sup_{t>0}\frac 1{t^r}\sup\ens{ \norm{x+ty}'+\norm{x-ty}'-2: \norm{x}'=\norm{y}'=1  }<\infty.
 \end{equation*}
\end{Definition}

From \cite{MR0407963}, we know that if $B$ is $r$-smooth and separable,
then there exists a constant $D$ such that for any sequence of
$B$-valued martingale differences $\prt{X_i}_{i\geq 1}$,

\begin{equation}\label{eq:rD_smooth_martingale}
  \mathbb{E}  \Big\|\sum_{i=1}^nX_i\Big\|^r  \leq D  \sum_{i=1}^n \mathbb{E} \norm{X_i}^r .
\end{equation}

\begin{Definition}
Let $1<r\leq 2$ and $D>0$.
We say that a Banach space $B$ is $\prt{r,D}$-smooth if $B$ is $r$-smooth and
inequality
\eqref{eq:rD_smooth_martingale} holds for all $B$-valued martingale difference
sequences $\prt{X_i}_{i\geq 1}$.
\end{Definition}

\section{Main results}  \label{sec2}

We start by a result for martingales difference  sequences
whose tail
have a uniform exponential decay.

\begin{theorem}
\label{th1} Let $\alpha \in (0, 1)$. Assume that $(X_{i}, \mathcal{F}_{i})_{i\geq 1}$ is a sequence of
 martingale differences in a $(2, D)$-smooth separable Banach space and  satisfies
 \begin{equation}
  \sup_{ i \geq 1}
 \sup_{t>0}\exp\left\{t^{\frac{2\alpha}{1-\alpha}}\right\}
 \PP\prt{\norm{X_i}>t}  \leq C_1
 \end{equation}
 for some constant $C_1$. Then, for all $ x>0$,
\begin{eqnarray}
\mathbb{P}\left( \max_{1\leq k \leq n}\| S_k \| >n x \right) &\leq& C(\alpha, x)  \exp\left\{-\left(\frac{ x  }{ 4D }\right)^{2\alpha} n^\alpha\ \right\} ,
\end{eqnarray}
where
\[
C(\alpha,  x)=  2+ 1007156 e^2D^2 C_1  \left( \frac{1}
{ x^{2\alpha} 16^{1-\alpha}D^{2(1-\alpha)}}  + \frac{1}{  x^2 }
\left( \frac{3(1-\alpha)}{2\alpha}\right)^{\frac{ 1-\alpha}{\alpha}} \right).
\]
does not depend  on $n$.
In particular, with $x=1$, it holds
\begin{eqnarray}
\mathbb{P}\left( \max_{1\leq k \leq n}\| S_k \| >n   \right)  &=&O \left( \exp\ens{ - \frac{ 1 }{ (4D)^{2\alpha} }  \, n^{\alpha} } \right),\ \ \ \ n\rightarrow \infty. \label{fnk}
\end{eqnarray}
\end{theorem}

\begin{remark}
 The condition on the decay on the tail can be rewriten
 as \begin{equation}
     \sup_{ i \geq 1}
 \sup_{s>0}s
 \PP\prt{ \exp\left\{ \norm{X_i}^{\frac{2\alpha}{1-\alpha}}\right\}>s}  \leq C_1.
    \end{equation}
In particular, it does not imply finiteness of
$\E{\exp\left\{ \norm{X_i}^{\frac{2\alpha}{1-\alpha}} \right\}  }$,
which was assumed in Theorem~2.1 of \cite{MR3005732}.
\end{remark}

\begin{remark}
 In \cite{MR2568267}, for $\mathbb{P}\left( \max_{1\leq k \leq n}\| S_k \| >n  \right) $, a rate of order $e^{-n}$ can be
 obtained under the condition that for some $\delta>0$,
 $\E{\exp\left\{\delta\abs{X_i} \right\} \mid \mathcal{F}_{i-1}}\leq K$ almost
 surely for
 some constant $K$.
 On one hand, the obtained rate is better than the one we
 obtained. On the other hand, there are situations where
 our results apply but not the one in \cite{MR2568267}.
 For example, let $\prt{Y_i}_{i\geq 0}$ be an independent
 sequence of random variables where for $i\geq 1$,
 $\PP\prt{Y_i=1}=\PP\prt{Y_i=-1}=1/2$ and $Y_0$ is a
 non-bounded random variable such that $\sup_{t>0}\exp\left\{t^{\frac{2\alpha}{1-\alpha}}\right\}
 \PP\prt{\abs{Y_0}>t} $ is finite for some $\alpha\in\prt{0,1}$. Letting $\mathcal{F}_i:=\sigma\prt{Y_j,0\leq j\leq i}$
 and $X_i:= Y_0Y_i$, then $\prt{X_i,\mathcal{F}_i}_{i\geq 1}$ is a martingale difference  sequence. Since
 $\abs{X_i}=\abs{Y_0}$, it follows that
 $\E{\exp\left\{\delta\abs{X_i}  \right\} \mid \mathcal{F}_{i-1}}
 =\exp\left\{\delta\abs{Y_0}  \right\} $, which is not bounded.

\end{remark}

We also investigate the case of the martingale differences having polynomial tail probability.

\begin{theorem}
\label{th2}
Let $B$ be an $\prt{r,D}$-smooth separable
Banach space where $1<r\leq 2$.
Let $p_2\geq p_1 >r$. Assume that $\prt{X_{i}, \mathcal{F}_{i}}_{i\geq 1}$ is a sequence of
$B$-valued martingale differences  and  satisfies
 $\ds\sup_{ i\geq 1 } \sup_{t>0}t^{p_1}\PP\prt{\norm{X_i}>t}\leq C_1 $ and
 $\ds\sup_{i\geq 1}\sup_{t>0}t^{p_2}\PP\prt{ \prt{\E{    \norm{X_i}^r\mid \F_{i-1}} }^{1/r}
 >t}\leq C_2 $
 for some constants $C_1$ and $C_2$. Then, for all $ x>0$,
\begin{multline}
\mathbb{P}\left( \max_{1\leq k \leq n}\| S_k \| >n x \right)\\  \leq   K_1(p_1,p_2,r,D)
C_1 x^{-p_1}\frac{1}{n^{p_1- 1}}   +  K_2(p_1,p_2,r,D)C_2 x^{-p_2}\frac{ 1}{n^{p_2-p_2/r} }  ,
 \label{eq:deviation_martingale_polynomial}
\end{multline}
where
\[
 K_1(p_1,p_2,r,D)=  \frac{2^{2p_2}}{2^{2p_2}-1}  2^{1-r} 2^{p_1+2p_1p_2/r}D^{p_1/r} \ \ \ \ \  \textrm{and}
\]
\begin{equation}
 K_2(p_1,p_2,r,D)=  \frac{2^{2p_2}}{2^{2p_2}-1}2^{1-r}2^{p_2+2p_2^2/r}D^{p_2/r}
    \prt{ \frac{p_2}{p_2-r}}^{p_2/r}
\end{equation}
do not depend  on $n$ or $x$.
In particular, with $x=1$, it holds
\begin{eqnarray}
\mathbb{P}\left( \max_{1\leq k \leq n}\| S_k \| >n   \right)  &=&O
\left( \frac{ 1}{n^{  \min\ens{ p_1-1,    p_2-p_2/r}   } } \right),\ \ \ \ n\rightarrow \infty.
\end{eqnarray}
\end{theorem}

\begin{remark}
When $B=\mathbb R$ and $p_1=p_2=p$, the rate for $\mathbb{P}\left( \max_{1\leq k \leq n}\| S_k \| >n   \right) $ is $n^{-p/2 }$. We thus recover
 the optimal convergence rate   of \cite{MR1856684}. Notice that  Theorem~3.6 (up
to the constants) is stated under weaker assumptions, since we do not need a finite moment of
order $p$ for $X_i$.
\end{remark}

\begin{remark}
Theorem~2.3 in \cite{1603.00432} gives a similar result as our Theorem~\ref{th2} in the case
 $r=2$. The condition
therein is in appearence different, because the condition of boundedness of the
 weak-$\mathbb{L}^p$-moments is replaced by the existence of random variables $X$ and $V$ such that
 for all $i$ and all $t$, $\PP\prt{\abs{X_i}>t}\leq \PP\prt{X>t}$ and
 $\PP\prt{\E{\norm{X_i}^r\mid\F_{i-1}}   >t  } \leq \PP\prt{V>t}$. However, if $X$ has a
 finite weak-$\mathbb{L}^p$-moment, then $C:=
 \sup_{ i\geq 1 } \sup_{t>0}t^{p_1}\PP\prt{\norm{X_i}>t}$
 is finite. Conversely, if $C$ is finite, then we can assume, by rescaling, that $C=1$; then take $X$ such that
 $\PP\prt{X>t}=\min\ens{1,t^{-p}}$.

 What the result of Theorem~\ref{th2} brings is the following. First, the case of $r$-smooth Banach
 spaces is considered here, whereas in Theorem~2.3 in \cite{1603.00432}, only the case of
 $2$-smooth Banach spaces is considered. Second, in our result, the constants are explicit.
\end{remark}

Theorem~\ref{th2} can also be used for sequences of independent centered random variables, which
are particular cases of martingale differences. Since the random variables $\E{    \norm{X_i}^r\mid \F_{i-1}}
$ are constant, one can apply Theorem~\ref{th2} for any $p_2$. In particular, we
can choose $p_2$ such that the decay in $n$ in the right hand side of
\eqref{eq:deviation_martingale_polynomial} is the same for both terms, namely,
$p_2=\prt{p_1-1}r/\prt{r-1}$.

\begin{corollary}
Let $B$ be an $\prt{r,D}$-smooth separable
Banach space where $1<r\leq 2$.
Let $p >r$. Assume that $\prt{X_{i} }_{i\geq 1}$ is an independent sequence of
$B$-valued random variables and  satisfies
 $\sup_{ i\geq 1 } \sup_{t>0}t^{p}\PP\prt{\norm{X_i}>t}\leq C  $
 for some constant  $C$. Then, for all $ x>0$,
\begin{equation}
\mathbb{P}\left( \max_{1\leq k \leq n}\| S_k \| >n x \right)  \leq   K(p,r,D) \bigg(x^{-p}
 C     +
 \sup_{i\geq 1}\prt{ \mathbb{E} \norm{X_i}^r }^{p/r}   x^{-\prt{p-1}r/\prt{r-1}} \bigg)\frac{ 1}{n^{p-1} } ,
 \label{eq:deviation_indep_polynomial}
\end{equation}
where
\begin{equation}
K(p,r,D)=  \frac{2^{2p_2}}{2^{2p_2}-1}  2^{1-r} 2^{p +2p p_2/r}D^{p /r}, \  \ \textrm{with} \ \  p_2=\prt{p -1}r/\prt{r-1},
\end{equation}
does not depend  on $n$ or $x$. In particular, with $x=1$, it holds
\begin{equation}
\mathbb{P}\left( \max_{1\leq k \leq n}\| S_k \| >n x \right) =
O\prt{ \frac 1{n^{p-1}}   }.
\end{equation}
\end{corollary}

\begin{remark}
It has been shown in Proposition~2.6 in \cite{MR1856684} that the power $p-1$ is
optimal, even for i.i.d. sequences with a finite moment of order $p$.
\end{remark}

\section{Proof of Theorems}

\subsection{Proof of Theorem~\ref{th1}}\label{subsec2}

The proof will be done by a truncation argument, similar method for
univariate martingale differences can be found  in Lesigne and Voln\'{y} \cite{MR1856684}.
For the bound part, we shall need the following Pinelis' inequality (cf.\ Theorem 3.5 of Pinelis \cite{MR1331198}).
\begin{lemma} \label{lemma1}
Assume that $\prt{X_{i}, \mathcal{F}_{i}}_{i\geq 1}$ is a sequence of
 martingale differences in a $(2, D)$-smooth separable Banach space and  satisfies
 $\norm{ X_{i}}_\infty \leq b$ for all $i\geq 1$.
 Then, for all $ x \geq 0$,
\begin{eqnarray}
\mathbb{P}\left( \max_{1\leq k \leq n} \| S_k \|  \geq x  \right)&\leq&
 \exp\left\{ - \frac{ x^2}{2D^2nb^2 }\right\} . \label{hoemax}
\end{eqnarray}
\end{lemma}

Let $\prt{X_{i}, \mathcal{F}_{i}}_{i \geq 1}$ be a sequence of martingale differences in  a $(2, D)$-smooth separable Banach space.   Given $u> 0$, define
\begin{eqnarray*}
X'_{i} &=& X_{i}\mathbf{1}_{\{\|X_{i}\|\leq u\}}-\mathbb{E}[X_{i}\mathbf{1}_{\{\|X_{i}\|\leq u\}}|\mathcal{F}_{i-1}],\\
X''_{i} &=& X_{i}\mathbf{1}_{\{\|X_{i}\|> u\}}-\mathbb{E}[X_{i}\mathbf{1}_{\{\|X_{i}\|> u\}}|\mathcal{F}_{i-1}], \\
S'_k &=&\sum_{i=1}^k X'_{i},  \ \ \ \ \ \ \ \ \ S''_k\ =\ \sum_{i=1}^k X''_{i}.
\end{eqnarray*}
Then $\prt{X'_{i}, \mathcal{F}_{i}}_{i\geq 1}$ and $(X''_{i}, \mathcal{F}_{i})_{i\geq 1}$ are two martingale difference sequences in  a $(2, D)$-smooth separable Banach space  and
$S_k= S'_k + S''_k $.
Let $t \in (0,1).$ For any $x>0$, it holds
\begin{eqnarray}
\mathbb{P}\left( \max_{1\leq k \leq n}\| S_k \| > x \right)
&\leq& \mathbb{P}\left( \max_{1\leq k \leq n}\| S'_k \| > x t \right) + \mathbb{P}\left( \max_{1\leq k \leq n}\|S''_k\|> x (1-t) \right). \label{fsum}
\end{eqnarray}
Using Lemma \ref{lemma1} and the fact that $\norm{X'_{i}}
\leq 2u$, we get
\begin{eqnarray}
\mathbb{P}\left( \max_{1\leq k \leq n} \| S'_k \| > x t \right)&\leq&  \exp\left\{-\frac{ x^2t^2}{8D^2 n u^2} \right\} . \label{f7i}
\end{eqnarray}
Using Theorem 4.1 of Pinelis \cite{MR1331198}, we get
\begin{eqnarray}
\mathbb{P}\left( \max_{1\leq k \leq n} \| S''_k \| > x (1-t) \right) &\leq& \frac{7200e^2 D^2 }{x^2(1-t)^2}\sum_{i=1}^n \mathbb{E} \norm{ X''_{i}}^{2}. \label{fk}
\end{eqnarray}
Let $F_{i}(x)=\mathbb{P}( \|X_{i} \| > x), x\geq 0.$
By assumption we have, for all $i$ and $x\geq0$,
\[
F_{i}(x)    \leq C_1 \exp\{ -x^{\frac{2\alpha}{1-\alpha}}\}.
\]
Using the inequality
\begin{equation}
 \E{\norm{X-\E{X\mid \mathcal G}}^2}\leq 4\E{\norm{X}^2},
\end{equation}
one gets that $\mathbb{E} \norm{ X''_{i}}^{2}
\leq 4\E{ \norm{ X_{i}}^{2}\mathbf{1}_{\{\|X_{i}\|> u\}}}$.
It is easy to see that
\begin{eqnarray}
 \mathbb{E}\|X'' _i \| ^2&\leq & - 4 \int_{u}^\infty  t^2dF_{i}(t)\nonumber\\
  &= & 4 u^2F_{i}(u) +  \int_{u}^\infty 8 tF_{i}(t) dt  \nonumber\\
&\leq & 4C_1 u^2\exp\{ -u^{\frac{2\alpha}{1-\alpha}}\}+8C_1  \int_{u}^\infty t \exp\{ -t^{\frac{2\alpha}{1-\alpha}}\} dt. \label{fkl}
\end{eqnarray}
Notice that the function $g(t)=t^3 \exp\{ -t^{\frac{2\alpha}{1-\alpha}}\}$ is decreasing in $[\beta, +\infty)$ and is increasing in $[0, \beta]$, where $\beta =\left(\frac{3(1-\alpha)}{2\alpha} \right)^{\frac{1-\alpha}{2\alpha}}$.
If $0<  u < \beta$, we have
\begin{eqnarray}
 \int_{u}^\infty t \exp\{ -t^{\frac{2\alpha}{1-\alpha}}\} dt &\leq &  \int_{u}^\beta t \exp\{ -t^{\frac{2\alpha}{1-\alpha}}\} dt   + \int_{\beta}^\infty t^{-2} t^3 \exp\{ -t^{\frac{2\alpha}{1-\alpha}}\} dt \nonumber \\
 &\leq&  \int_{u}^\beta t \exp\{ -u^{\frac{2\alpha}{1-\alpha}}\} dt  + \int_{\beta}^\infty t^{-2} \beta^3 \exp\{ -\beta^{\frac{2\alpha}{1-\alpha}}\} dt \nonumber \\
 &\leq& \frac{3}{2}\beta^2\exp\{ -u^{\frac{2\alpha}{1-\alpha}}\}. \label{ghds}
\end{eqnarray}
If $ \beta \leq  u$, we have
\begin{eqnarray}
 \int_{u}^\infty t \exp\{ -t^{\frac{2\alpha}{1-\alpha}}\} dt &=&  \int_{u}^\infty t^{-2} t^3 \exp\{ -t^{\frac{2\alpha}{1-\alpha}}\} dt \nonumber \\
 &\leq&  \int_{u}^\infty t^{-2} u^3 \exp\{ -u^{\frac{2\alpha}{1-\alpha}}\} dt \nonumber \\
 &=& u^2\exp\{ -u^{\frac{2\alpha}{1-\alpha}}\}. \label{xvcb}
\end{eqnarray}
By  \eqref{fkl},  \eqref{ghds}  and \eqref{xvcb}, it follows that
\begin{eqnarray}
 \mathbb{E}\|X''_i\|^2 &\leq &12 C_1( u^2 +  \beta^2 ) \exp\ens{ -u^{\frac{2\alpha}{1-\alpha}}}.
\end{eqnarray}
From  \eqref{fk}, we get
\begin{eqnarray}
\mathbb{P}\left( \max_{1\leq k \leq n} \| S''_k \| > x (1-t) \right)&\leq& \frac{12\times 7200 e^2  D^2 C_1 n  }{x^2(1-t)^2} (u^2+ \beta^2)\exp\ens{ -u^{\frac{2\alpha}{1-\alpha}}}. \label{fl}
\end{eqnarray}
Combining \eqref{fsum}, \eqref{f7i} and \eqref{fl} together, we obtain
\begin{eqnarray}
\mathbb{P}\left( \max_{1\leq k \leq n} \| S_k \| >x \right) &\leq&   2\exp\left\{-\frac{ x^2t^2}{8D^2 u^2 n } \right\} +\frac{12\times 7200e^2 D^2 C_1 n  }{ (1-t)^2} \left( \frac{u^2}{x^2}+ \frac{\beta^2}{x^2} \right) \exp \ens{ -u^{\frac{2\alpha}{1-\alpha}}  }.\nonumber
\end{eqnarray}

Taking $t=\frac{1}{\sqrt{2}}$ and $u=  \left(\frac{x }{4D\sqrt{n}} \right)^{1-\alpha},$ we get, for all $x>0$,
\begin{eqnarray*}
\mathbb{P}\left( \max_{1\leq k \leq n} \| S_k \| >x \right) &\leq& C_n(\alpha,  x)\exp\left\{-\left(\frac{ x^2 }{ 16 D^2 n }\right)^\alpha\ \right\}, \nonumber
\end{eqnarray*}
where
\[
C_n(\alpha,   x) =  2+ \frac{86400}{\prt{1-1/\sqrt 2}^2} e^2 D^2 C_1 n \left( \frac{1}{x^{2\alpha}(16D^2 n)^{1-\alpha}}  + \frac{\beta^2}{x^2} \right) .
\]
Hence (using the fact that $\frac{86400}{\prt{1-1/\sqrt 2}^2}
\leq 1007156$), for all $x>0$,
\begin{eqnarray*}
\mathbb{P}\left( \max_{1\leq k \leq n} \| S_k\| >n x \right) &\leq&C(\alpha, x)  \exp\left\{-\left(\frac{ x  }{ 4D }\right)^{2\alpha} n^\alpha\ \right\}, \nonumber
\end{eqnarray*}
where
\[
C(\alpha,  x)=  2+1007156 e^2D^2 C_1 \left( \frac{1}
{ x^{2\alpha} 16^{1-\alpha}D^{2(1-\alpha)}}  + \frac{1}{  x^2 }
\left( \frac{3(1-\alpha)}{2\alpha}\right)^{\frac{ 1-\alpha}{\alpha}} \right).
\]
This completes the proof  of Theorem \ref{th1}.

\subsection{Proof of Theorem~\ref{th2}}

Theorem~\ref{th2} will be proven by an application of Theorem~1.3 in \cite{1603.00432},
which states the following.
\begin{lemma}
   Let $\prt{B,\norm{\cdot}}$ be a separable $\prt{r,D}$-smooth Banach space, where
   $1<r\leqslant 2$.
   For any $B$-valued martingale
   difference sequences $\prt{X_i,\F_i}_{i\geq 1}$, the following inequality
   holds for each $n\geq 1$ and $q, x>0$:
   \begin{equation}\label{eq:deviation_inequality_non_stationary_martingale_diff}
    \PP\prt{\max_{1\leqslant i\leqslant n}\norm{S_i}>x }
    \leq I_{1}(q, x)+I_{2}(q, x),
   \end{equation}
   where
   $$I_{1}(q, x)= \frac{2^{q-r}q }{2^q-1}\int_0^1\PP\prt{\max_{1\leqslant i\leqslant n}
    \norm{X_i}>2^{-1-q/r}D^{-1/r}xu}u^{q-1}\mathrm du,$$
   $$I_{2}(q, x)= \frac{2^{q-r}q}{2^q-1} \int_0^1\PP\bigg(
    \Big(\sum_{i=1}^n
    \E{\norm{X_i}^{r}\mid \F_{i-1}}\Big)^{1/r}>2^{-1-q/r}D^{-1/r}xu \bigg)u^{q-1}\mathrm du$$
    and $D$ is a constant satisfying
   \eqref{eq:rD_smooth_martingale} for any $n$ and any martingale difference sequences.
\end{lemma}

Applying inequality \eqref{eq:deviation_inequality_non_stationary_martingale_diff}  with $q=2p_2$ and $x$ replaced by $nx$, we get for each $n\geq 1$ and $p_2, x>0$:
\begin{equation}\label{eq:application_deviation_inequality_polynomial}
 \PP\prt{\max_{1\leqslant i\leqslant n}\norm{S_i}>nx }
    \leq   I_{1}(2p_2, nx)+I_{2}(2p_2, nx),
   \end{equation}
In order to control the terms appearing in the right hand side of \eqref{eq:application_deviation_inequality_polynomial}, we
introduce the following quantity, for $p>1$
\begin{equation}
N_p\prt{X}:=\sup_{t>0} t^p\PP\prt{\abs{X}>t}  .
\end{equation}
Then by definition, for all random variable $X$,
$\PP\prt{\abs{X}>t}\leq t^{-p}N_p\prt{X}$. Therefore,
\begin{align}
\PP\prt{\max_{1\leqslant i\leqslant n}
    \norm{X_i}>2^{-1-2p_2/r}D^{-1/r}nxu}
    &\leqslant \sum_{i=1}^n\PP\prt{
    \norm{X_i}>2^{-1-2p_2/r}D^{-1/r}nxu} \nonumber\\
   & \leq \sum_{i=1}^n\prt{2^{-1-2p_2/r}D^{-1/r}nxu}^{-p_1}N_{p_1}\prt{X_i}\nonumber\\
   & \leq  C_1 n^{1-p_1} 2^{p_1+2p_1p_2/r}D^{p_1/r} x^{-p_1}u^{-p_1}.
\end{align}
Using the last bound, we derive that
\begin{equation}\label{eq:proof_polynomial_first_term}
 I_{1}(2p_2, nx) \leq  \frac{2^{2p_2}}{2^{2p_2}-1}  2^{1-r} 2^{p_1+2p_1p_2/r}D^{p_1/r} C_1 n^{1-p_1} x^{-p_1}.
\end{equation}
In order to control the second term of \eqref{eq:application_deviation_inequality_polynomial}, 
 we would like to use the uniform control on
$N_{p_2}\prt{\prt{  \E{ \norm{X_i}^r \mid \F_{i-1}   }   }^{1/r}  }$. However, the
triangle inequality fails for $N_{p_2}$. This leads us to introduce the weak-$\mathbb L^p$-norm
defined for $p>1$ as
\begin{equation}
\norm{X}_{p,\infty}:= \sup\ens{\PP\prt{A}^{-1+1/p}\E{\norm{X}\mathbf 1_A},
A\in \F, \PP\prt{A}>0}.
\end{equation}
Then $\norm{\cdot}_{p,\infty}$ defines a norm and is linked to $N_p$ in the following way:
\begin{equation}\label{eq:inequality_weak_Lp}
N_p\prt{X}\leq \norm{X}_{p,\infty}\leq \frac p{p-1}N_p\prt{X}.
\end{equation}
The first inequality follows from the estimate
\begin{equation}
t^p \PP\ens{\norm{X}>t}\leq
\E{\norm{X}\mathbf 1\ens{\norm{X}>t}  }\leq
\norm{X}_{p,\infty}\prt{\PP\ens{\norm{X}>t}}^{1-1/p}.
\end{equation}
For the second one, fix $A\in\F$ such that $\PP\prt{A}>0$.
We write
\begin{equation}
 \E{\norm{X}\mathbf 1_A}=
 \int_0^{+\infty}\PP\prt{\ens{\norm{X}>t    }\cap A}dt
 \leq \int_0^{+\infty}\min\ens{\PP\prt{\norm{X}>t }, \PP\prt{A}
 }dt,
\end{equation}
bound $\PP\prt{\norm{X}>t }$ by $t^{-p}
N_p\prt{X}^p$ and compute the remaining integral.

It follows from \eqref{eq:inequality_weak_Lp} that
for a fixed $y>0$,
\begin{align*}
\PP\prt{
    \prt{\sum_{i=1}^n
    \E{\norm{X_i}^{r}\mid \F_{i-1}}}^{1/r}>y}& =
    \PP\prt{
    \sum_{i=1}^n
    \E{\norm{X_i}^{r}\mid \F_{i-1}} >y^r}\\
    &\leq y^{-p_2}\prt{ N_{p_2/r}
    \prt{ \sum_{i=1}^n
    \E{\norm{X_i}^{r}\mid \F_{i-1}}      }}^{p_2/r}\\
    &\leq y^{-p_2}
    \norm{ \sum_{i=1}^n
    \E{\norm{X_i}^{r}\mid \F_{i-1}}      }_{p_2/r,\infty}^{p_2/r}\\
    &\leq y^{-p_2}  \prt{
 \sum_{i=1}^n \norm{
    \E{\norm{X_i}^{r}\mid \F_{i-1}}      }_{p_2/r,\infty}}^{p_2/r}\\
    &\leq y^{-p_2} \prt{ \frac{p_2}{p_2-r}}^{p_2/r}
 \prt{\sum_{i=1}^n N_{p_2/r}\prt{
    \E{\norm{X_i}^{r}\mid \F_{i-1}}      }}^{p_2/r} ,
\end{align*}
which leads to
\begin{equation}
\PP\prt{
    \prt{\sum_{i=1}^n
    \E{\norm{X_i}^{r}\mid \F_{i-1}}}^{1/r}>y}\leq
    y^{-p_2}C_2\prt{ \frac{p_2}{p_2-r}}^{p_2/r}n^{p_2/r}.
\end{equation}
Applying this to $I_{2}(2p_2, nx)$ with $y= 2^{-1-2p_2/r}D^{-1/r}xun$ and integrating in $u$ gives
\begin{equation}\label{eq:proof_polynomial_second_term}
I_{2}(2p_2, nx) \leq \frac{2^{2p_2}}{2^{2p_2}-1}2^{1-r}2^{p_2+2p_2^2/r}D^{p_2/r}
    \prt{ \frac{p_2}{p_2-r}}^{p_2/r} C_2
    x^{-p_2}n^{p_2/r-p_2}.
\end{equation}
The combination of \eqref{eq:application_deviation_inequality_polynomial},
\eqref{eq:proof_polynomial_first_term} and \eqref{eq:proof_polynomial_second_term}
completes the proof of Theorem~\ref{th2}.

\def\polhk\#1{\setbox0=\hbox{\#1}{{\o}oalign{\hidewidth
  \lower1.5ex\hbox{`}\hidewidth\crcr\unhbox0}}}\def\cprime{$'$}
  \def\polhk#1{\setbox0=\hbox{#1}{\ooalign{\hidewidth
  \lower1.5ex\hbox{`}\hidewidth\crcr\unhbox0}}} \def\cprime{$'$}
\providecommand{\bysame}{\leavevmode\hbox to3em{\hrulefill}\thinspace}
\providecommand{\MR}{\relax\ifhmode\unskip\space\fi MR }
\providecommand{\MRhref}[2]{%
  \href{http://www.ams.org/mathscinet-getitem?mr=#1}{#2}
}
\providecommand{\href}[2]{#2}

\end{document}